\documentclass[12pt]{article}
\usepackage{fancyhdr, array,calc,graphicx,url,tabularx}
\usepackage{geometry}
\usepackage{latexsym}
\usepackage{amssymb}
\usepackage{amsmath}
\usepackage{makeidx}
\usepackage{enumerate}
\usepackage{verbatim}
\usepackage{tikz}
\usepackage[colorlinks=true,linkcolor=blue]{hyperref}

\usepackage{changepage}

\usepackage[utf8]{inputenc}
\usepackage{tikz}
\tikzset{
     block/.style={rectangle, draw, fill=red!40, text width=6em,
                   text centered, rounded corners, minimum height=3em},
     arrow/.style={-{Stealth[]}}
     }

\makeindex

{
   \end{minipage}
   \vspace*{\stretch{3}}
   \clearpage
}

\def\undertilde#1{{\baselineskip=0pt\vtop
  {\hbox{$#1$}\hbox{$\scriptscriptstyle\sim$}}}{}}

\newcommand{\utilde}{\undertilde}

\renewcommand{\gg}{\gamma}

\newcommand{\bR}{{\mathbb{R}}}

\newcommand{\card}[1]{{\vert #1 \vert} }

\renewcommand{\models}{\vDash}
\newcommand{\powerset}{{\wp}}
\newtheorem{theorem}{Theorem}[section]
\newtheorem{proposition}[theorem]{Proposition}

\newtheorem{lemma}[theorem]{Lemma}
\newtheorem{corollary}[theorem]{Corollary}

\newtheorem{conjecture}[theorem]{Conjecture}

\numberwithin{figure}{section}

\newenvironment{proof}{{\it{
Proof.}}}{\nopagebreak\mbox{}{\hfill$\square$}
\par\bigskip}
\newenvironment{Claim}{\noindent{\it
Claim} {}}{}
\newenvironment{Case}{\noindent{\it
Case}}{}

\newcommand{\rcon}[1]{Conjecture~\ref{#1}}

\newcommand{\rprop}[1]{Proposition~\ref{#1}}

\newcommand{\rlem}[1]{Lemma~\ref{#1}}

\newcommand{\rcor}[1]{Corollary~\ref{#1}}

\def\k{\kappa}
\def\a{\alpha}
\def\b{\beta}
\def\d{\delta}

\def\l{\lambda}

\def\Q{{\mathcal{ Q}}}

\def\R{{\mathcal R}}

\def\H{{\rm{HOD}}}

\def\S{{\mathcal{S}}}

\def\card#1{\left|#1\right|}

\def\and{\mathrel{\kern1pt\&\kern1pt}}

\def\<#1>{\langle\,#1\,\rangle}

 \input xy
 \xyoption{all}
       
\title{An inner model theoretic proof of Becker's theorem}

\author{Grigor Sargsyan}

\date{\today}

\begin{document}
\maketitle

\begin{abstract} We re-prove Becker's theorem from \cite{Becker} by showing that $AD^{L(\bR)}$ implies that $L(\bR)\models ``\omega_2$ is $\utilde{\delta}^2_1$-supercompact". Our proof uses inner model theoretic tools instead of Baire category. We also show that $\omega_2$ is $<\Theta$-strongly compact.
\end{abstract}

This article draws inspiration from the work of Neeman (\cite{IMU}) who, using inner model theoretic tools, showed that under $AD^{L(\bR)}$, $\omega_1$ is $<\Theta$-supercompact. We have also been influenced by the work of Becker (\cite{Becker}), Becker-Jackson (\cite{BeckerJackson}) and Jackson (\cite{JacksonSquare}). In \cite{Becker}, Becker showed that assuming $AD+V=L(\bR)$, $\omega_2$ is $\delta^2_1$-supercompact. In \cite{BeckerJackson}, Becker and Jackson showed that, under $AD+V=L(\bR)$, all projective cardinals are $\delta^2_1$-supercompact. Finally, in \cite{JacksonSquare}, Jackson showed that under $AD+V=L(\bR)$ all Suslin cardinals and their successors are $\delta^2_1$-supercompact.  

In this short note, we re-prove Becker's theorem using inner model theoretic tools. The paper assumes familiarity with what is commonly called $\H$ analysis. The reader can find this background exposited in \cite{IMU} and in \cite{SteelWoodin}. The point of re-proving such results is to find more applications of inner model theory in descriptive set theory. In particular, we strongly believe that connecting iteration sets with Kechris-Woodin generic codes will yield many applications, and thus invite the community to consider \rcon{conjecture}.

The author would like to thank the referee for noticing many typos and suggesting important improvements. The author's work was partially supported by the NSF Career Award DMS-1352034.

\section{Measures on $\powerset_{\omega_2}(\l)$}

We do not want to make the paper artificially long. The paper is aimed at experts of inner model theory, those who are familiar with the terminology of \cite{SteelWoodin}. 

We assume $AD+V=L(\bR)$. Fix $\l<\Theta$. Let $A$ be an $OD$ set of reals such that $\gg_{A, \infty}\geq\l$. Suppose $\R$ is a suitable premouse that is $A$-iterable. It is customary to let $\d^\R$ be the Woodin cardinal of $\R$. 
Assume that $\l\in rng(\pi_{(\R, A), \infty})$. We then let $\l^\R$ be such that $\pi_{(\R, A), \infty}(\l^\R)=\l$. 

We let $Code(A, \l)\subseteq \bR$ be the set of reals $x$ such that $x$ codes a pair $(\R_x, \a_x)$ such that $\R_x$ is an $A$-iterable suitable pre-mouse such that $\l^\R$ is defined  and $\a_x<\l^\R$. Let $\leq_{A, \l}$ be the natural pre-wellordering of $Code( A, \l)$ given by: $x\leq_{ A, \l} y$ if and only if whenever $S$ is an $A$-iterate of $\R_x$ and an $A$-iterate of $\R_y$, $\pi_{(\R_x, A), (\S, A)}(\a_x)\leq \pi_{(\R_y, A), (\S, A)}(\a_y)$. We have that $\leq_{ A, \l}$ has length $\l$. 
Given $x\in Code(A, \l)$ let
\begin{center}
 $c(x)=\pi_{(\R_x, A), \infty}(\a_x)=\card{x}_{\leq_{ A, \l}}$.
 \end{center}

Let $S$ be a tree of a $\Sigma^2_1$-scale on a universal $\Sigma^2_1$-set. Given $x$ and $y$ we write $x\sim_S y$ if and only $x\in L[S, y]$ and $y\in L[S, x]$. We then say that $d$ is an $S$-degree if $d$ is an $\sim_S$-class. We write $d\leq_S e$ if  $d\in L[S, e]$. Let now $C(A, \l)=\{ d: Code(A, \l)\cap HC^{L[S, d]}\not =\emptyset\}$. The following are two key points to keep in mind:
\begin{enumerate}
\item $\sim_S$ is an equivalence relation,
\item $C(A, \l)$ contains an $S$-cone, i.e., there is an $S$-degree $e$ such that whenever $e\leq_S d$, $d\in C(A, \l)$.
\end{enumerate}

The following is a corollary to the Harrington-Kechris theorem (see \cite{HarKech}, and see \cite{ADRUB} and the references there for some uses of it).

%


 
 \begin{corollary}\label{cor 1} There is a formula $\phi$ such that whenever $d\in C(A, \l)$, $g$ is $<\omega_1^V$-generic over $L[S, d]$ and $\R\in L_{\omega_1^V}[S, d][g]$, 
 \begin{center}
 $\R$ is a suitable pre-mouse if and only if $L[S, d][g]\models \phi[\R]$.  
 \end{center}
 Moreover, there is a formula $\psi$ such that for any $A$-iterable suitable $\Q, \R\in L_{\omega_1^V}[S][g]$ and for any $\pi$, 
  \begin{center}
 $\R$ is an $A$-iterate of $\Q$ and $\pi: H^\Q_A\rightarrow H^\R_A$ is the $A$-iteration embedding if and only if $L[S, d][g]\models \psi[\Q, \R, \pi, \tau_A]$,  
 \end{center}
 where $\tau_A$ is the term relation for $A$ in $L[S, d]^{Coll(\omega, <\omega_1^V)}$.
 \end{corollary}
 The formulas $\phi$ and $\psi$ essentially repeat the definitions of suitability and $A$-iterability. Another important lemma that we need is a consequence of what is usually called \textit{generic comparisons} (see \cite{SteelWoodin}). The proof is a standard generic comparison argument which we leave to the reader.
 
 \begin{lemma}\label{generic comparisons} Suppose $d\in C(A, \l)$ and $g$ is $<\omega_1^V$-generic over $L[S, d]$. Suppose $\phi$ is as in \rcor{cor 1}, and for some $\R\in L[S, d][g]$, $L[S, d][g]\models \phi[\Q, \R]$. Then there is an $\emptyset$-iterate $\S$ of $\R$ such that $\S\in L_{\omega_1^V}[S, d]$.  
 \end{lemma}
 
 Given $d\in C(A, \l)$ we let $B_d$ be the set of $\b$ such that there is $x\in Code(A, \l)$ with the property that $(\R_x, \a_x)\in L[S, d]$ and $c(\R_x, \a_x)=\b$. As $\card{L_{\omega^V_1}[S, d]}=\omega^V_1$, we have that $B_d\in \powerset_{\omega_2}(\l)$.  \rlem{generic comparisons} has the following easy corollary.

\begin{corollary}\label{cor 2} Suppose $d_0\in C(A, \l)$ and $d$ is a $S$-degree such that $L[S, d]$ is a $<\omega_1^V$-generic extension of $L[S, d_0]$. Then $B_{d_0}=B_d$.
\end{corollary}

We now define $\mu(A, \l)$ on $\powerset_{\omega_2}(\l)$ by setting $B\in \mu(A, \l)$ if and only if for an $S$-cone of $d$, $B_d\in B$. 

\begin{lemma}\label{completeness} $\mu(A, \l)$ is an $\omega_2$-complete ultrafilter on $\powerset_{\omega_2}(\l)$. 
\end{lemma}
\begin{proof} Clearly $\mu(A,\l)$ is an ultrafilter. Let $(B_\xi: \xi<\omega_1)$ be such that $B_\xi\in \mu(A, \l)$ for all $\xi<\omega_1$. Let $WO$ be the set of reals coding a countable ordinal. Using the coding lemma we can find $y\in \bR$ and a $\Sigma^1_2(y)$-set $D\subseteq WO\times \bR$ such that 
\begin{enumerate}
\item $[y]_S\in C(A, \l)$,
\item for every $x\in WO$, $D_x\not=\emptyset$ (here $D_x=\{ z: (x, z)\in D\}$),
\item for every $x\in WO$, $D_x\subseteq \{z: [z]_S$ is a base of a cone witnessing that $B_{\card{x}}\in \mu(A, \l)\}$\footnote{where $[z]_S$ is the $S$-degree given by $z$ and $\card{x}$ is the ordinal coded by $x$.}.
\end{enumerate}
Let $d\in C$ be such that $y\in L[S, d]$. We claim that $B_d\in B_\xi$ for every $\xi<\omega_1$. To see this, fix $\xi<\omega_1$. Let $g\subseteq Coll(\omega, \xi)$ be $L[S, d]$-generic and $u$ be a real such that $L[S, d][g]=L[S, u]$. Let $x\in \bR^{L[S, u]}$ be such that $\card{x}=\xi$. Because $D$ is $\Sigma^1_2(y)$ we have that there is $z\in D_x\cap L[S, u]$. Because $[z]_S\leq [u]_S$, we must have that $B_{[u]_S}\in B_\xi$. However, it follows from \rcor{cor 2} that $B_d=B_{[u]_S}$. Hence, $B_d\in B_\xi$. 

As $d$ was arbitrary, we have shown that for any $d$ that is $S$-above $[y]_S$, $B_d\in \cap_{\xi<\omega_1}B_\xi$. It follows that $\cap_{\xi<\omega_1}B_\xi\in \mu(A, \l)$.
\end{proof}

 \section{$\omega_2$ is $\utilde{\delta}^2_1$-supercompact and $<\Theta$-strongly compact}
 
 \begin{proposition}\label{fines} For every $\l<\Theta$ and an ordinal definable $A\subseteq \bR$ such that $\gg_{A, \infty}\geq \l$, $\mu(A, \l)$ is superfine, i.e.,
 for every $B\in \powerset_{\omega_2}(\l)$, 
 \begin{center}
 $\{ D\in \powerset_{\omega_2}(\l): B\subseteq D\}\in \mu(A, \l)$.
 \end{center}
 \end{proposition}
 \begin{proof} Fix $B$ and let $f:\omega_1\rightarrow B$ be a bijection. Let $B_\xi=\{ x\in Code(A, \l): c(x)=f(\xi)\}$. Using the coding lemma find $y\in \bR$ and $D\subseteq WO\times \bR$ such that
 \begin{enumerate}
 \item $[y]_S\in C(A, \l)$,
\item $D\in \Sigma^1_2(y)$,
 \item for every $x\in WO$, $D_x\not =\emptyset$,
 \item for every $x\in WO$, $D_x\subseteq \{z\in Code(A, \l): c(z)=f(\card{x})\}$.
 \end{enumerate}
 We claim that for every $d$ such that $[y]_S\leq_S d$, $B\subseteq B_d$. To see this, fix $d$ such that $[y]_S\leq_S d$. Fix $\zeta\in B$. We want to see that $\zeta\in B_d$. Let $\xi=f^{-1}(\zeta)$, and fix $u\in \bR$ such that $L[S, u]$ is a generic extension of $L[S, d]$ and $\xi$ is countable in $L[S, u]$. Fix $x\in WO\cap L[S, u]$ such that $\card{x}=\xi$. Because $D\in \Sigma^1_2(y)$, we have that $D_x\cap L[S, u]\not=\emptyset$. Fix then $z\in Code(A, \l)\cap D_x\in L[S, u]$. It follows that $c(z)=f(\xi)$. Since $c(z)\in B_{[u]_T}=B_d$, we have that $\zeta\in B_d$. 
 \end{proof}
 
 Putting \rprop{completeness} and \rprop{fines} we get the following corollary.

\begin{corollary} Assume $AD+V=L(\bR)$. Then $\omega_2$ is a $<\Theta$-strongly compact. More precisely, for every $\l<\Theta$ there an $\omega_2$-complete superfine ordinal definable ultrafilter on $\powerset_{\omega_2}(\l)$. 
\end{corollary} 

\begin{theorem}[Becker, \cite{Becker}] Assume $AD+V=L(\bR)$. Then $\omega_2$ is $\utilde{\delta}^2_1$-supercompact.
\end{theorem}
\begin{proof}
Set $\l=\utilde{\delta}^2_1$. Suppose $\R$ is an $\emptyset$-iterable suitable pre-mouse. Recall that if $\nu$ is the least cardinal that is $<\d^\R$-strong in $\R$ then $\pi_{(\R, \emptyset), \infty}(\nu)=\l$ (see \cite[Chapter 8]{OIMT}).  We now want to show that $\mu=_{def}\mu(\emptyset, \l)$ is an $\omega_2$-supercompactness measure. \rprop{completeness} shows that $\mu$ is $\omega_2$-complete and \rprop{fines} shows that $\mu$ is fine. It remains to show that $\mu$ is normal. The following lemma is the first step towards normality. Set $Code=_{def}Code(\emptyset, \l)$ and $\leq^*=\leq_{\emptyset, \l}$.

\begin{lemma}\label{descending function} Suppose $F:\powerset_{\omega_2}(\l)\rightarrow \l$ is such that for an $S$-cone of $d$, $F(B_d)\in B_d$. Then for an $S$-cone of $d$ there is $x\in (\bR^{L[S, d]}\cap Code)$ such that $c(x)=F(B_d)$. 
\end{lemma}
\begin{proof} Assume not. Fix an $S$-degree $d_0$ such that whenever $d$ is $S$-above $d_0$, for every $x\in (\bR^{L[S, d]}\cap Code)$,  $c(x)\not =F(B_d)$. Fix $(\R, \a)\in L[S, d_0]$ such that $\pi_{(\R, \emptyset), \infty}(\a)=F(B_{d_0})$. 

Let $\nu<\omega_1$ be any cardinal of $L[S, d_0]$ such that $(\R, \a)\in L_\nu[S, d_0]$ and let $g\subseteq Coll(\omega, (\nu^+)^{L[S, d_0]})$ be $L[S, d_0]$-generic. Let $x\in \bR$ be such that $L[S, d_0][g]=L[S, x]$. We then have that $B_{d_0}=B_{[x]_S}$ (see \rcor{cor 1}). This is a contradiction as we can find $y\in L[S, x]\cap \bR$ coding $(\R, \a)$. 
\end{proof}

\begin{lemma} $\mu$ is normal. 
\end{lemma}
\begin{proof}
Suppose $\mu$ is not normal. Let $F:\powerset_{\omega_2}(\l)\rightarrow \l$ be such that for an $S$-cone of $d$, $F(B_d)\in B_d$ but $F$ is not constant on a $\mu$-measure one set. Let $e_0\in C$ be a base for the cone of the previous sentence. 

 Let $e\in C$ be $S$-above $e_0$ and such that for every $d$ such that $e\in L[S, d]$, there is $x\in (\bR^{L[S, d]}\cap Code)$ with the property that $c(x)=F(B_d)$.  We now follow an idea of Becker from \cite{Becker}.

Given an ordinal $\xi<\l$ let $D_\xi=\{ d: F(B_d)\not =\xi\}$. We have that for each $\xi$, $D_\xi$ contains an $S$-cone. Let then $C_\xi=\{ x\in \bR: [x]_S$ is a base of a cone contained in $D_\xi\}$. It follows from the coding lemma that there is a real $y$ and a set $D$ such that
\begin{enumerate}
\item $e\leq_S[y]_S$,
\item $H\subseteq Code\times \bR$ is $\Sigma^2_1(y)$,
\item if $(x, z)\in H$ then $z\in C_{c(x)}$,
\item for every $x\in Code$ there is $z$ such that $(x, z)\in H$. 
\end{enumerate}

Set $d=[y]_S$. Let $x\in Code\cap L[S, d]$ be such that $c(x)=F(B_d)$. Because $H$ is $\Sigma^2_1(y)$ we have that $H_x\cap L[S, d]\not=\emptyset$. Let then $z\in H_x\cap L[S,d]$. It follows that $[z]_S\leq_S d$. Hence, $d\in H_{c(x)}$. It follows that $F(B_d)\not =c(x)=F(B_d)$, contradiction. 
\end{proof}
\end{proof}

\section{A covering conjecture}

Again we assume $AD+V=L(\bR)$. Suppose $\k<\l<\Theta$ and $A$ is an ordinal definable set of reals such that $\gg_{A, \infty}\geq \l$. Given $X\in \powerset_{\k}(\l)$ we say that $X$ is an $A$-\textit{iteration set} if for every $\a \in X$ there is an $A$-iterable $\Q$ such that $\a \in rng(\pi_{(\Q, A), \infty})$ and $\pi_{(\Q, A), \infty}[\gg^\Q_A]\cap \l\subseteq X$.

\begin{conjecture}\label{conjecture} Assume $AD+V=L(\bR)$. Suppose $\kappa$ is either a Suslin cardinal or a successor of a Suslin cardinal. Then for every $\l<\Theta$, an OD set $A\subseteq \bR$ such that $\l\leq \gg_{A, \infty}$, and $B\in \powerset_\k(\l)$ there is an $A$-iteration set $X\in \powerset_\k(\l)$ such that $B\subseteq X$.
\end{conjecture}

The conjecture is clearly true for $\k=\omega_1$. \rprop{fines} shows that the conjecture is true for $\k=\omega_2$. We expect that the validity of the full conjecture will follow once a link is made between Kechris-Woodin generic codes and iteration sets.

\bibliographystyle{plain}
\bibliography{Omega_2-sc.bib}
\end{document}